\newcommand{\briefrefs}{T}

\if T\briefrefs

\else

\fi

\newcommand{\beqn}{\begin{equation}}
\newcommand{\eeqn}{\end{equation}}

\newcommand{\be}{\begin{equation}}
\newcommand{\ee}{\end{equation}}
\newcommand{\beqa}{\begin{eqnarray}}
\newcommand{\eeqa}{\end{eqnarray}}
\newenvironment{vectra}[1]%
	{\left[ 		
	\begin{array}{#1}}
	{\end{array} \right]}

\newcommand{\bv}{\begin{vectra}}
\newcommand{\ev}{\end{vectra}}

\newcommand{\ba}{\begin{array}}
\newcommand{\ea}{\end{array}}

\newcommand{\define}{\stackrel{\triangle}{=}}

\newcommand{\Comment}[1]{}

\newcommand{\figureLongCapTwoCol}[3]{
	\if F\draft
	\begin{figure}[htbp]
	\centerline{\psfig{figure=#2,width=3.25in}} 
	\caption{\protect {\small #3}}
	\label{#1} 
	\end{figure}
	\fi
}

\newcommand{\draft}{F}

\if T\draft
\documentclass[11pt,twoside,draft]{IEEEtran}
\usepackage{epsf}
\else
\documentclass[twocolumn]{IEEEtran}
\usepackage{epsf}
\fi

\begin{document}
\if T\draft
\bibliographystyle{IEEE}
\else
\bibliographystyle{IEEE}
\fi

\author{Federico Fontana and Davide Rocchesso\thanks{F.~Fontana is with the Dipartimento Scientifico e Tecnologico, Universit\`a degli Studi di Verona, 37134 Verona, Italy, (e-mail: fontana@sci.univr.it). D.~Rocchesso is with the Dipartimento Scientifico e Tecnologico,
Universit\`a degli Studi di Verona, 37134 Verona, Italy, 
(e-mail: rocchesso@sci.univr.it; http://www.sci.univr.it/\~{}rocchess) 
~~~{\copyright 2000 IEEE. Personal use of this material is permitted. However, permission to
                          reprint/republish this material for advertising or promotional purposes or for creating new
                          collective works for resale or redistribution to servers or lists, or to reuse any copyrighted
                          component of this work in other works must be obtained from the IEEE.}
}}

\markboth{IEEE Trans. on Speech and Audio, vol. 9, no. 2, february 2001.}{Fontana and Rocchesso: Waveguide Mesh Geometries}

\title{Signal-Theoretic Characterization of Waveguide Mesh Geometries for Models of Two--Dimensional Wave Propagation in Elastic Media}

\date{\today}
\Comment{
{
\begin{center}
{\Large Signal-Theoretic Characterization of Waveguide Mesh Geometries for Models of Two--Dimensional Wave Propagation in Elastic Media}\\
{Federico Fontana\footnote{F.~Fontana is with the Centro di Sonologia Computazionale, Universit\`a degli Studi di Padova, via S.~Francesco 11, 35121 Padova - ITALY, E-mail fefo@csc1.unipd.it.  } and Davide Rocchesso\footnote{D.~Rocchesso (corresponding author) is with the Dipartimento Scientifico e Tecnologico,
Universit\`a degli Studi di Verona, strada Le Grazie, 37134 Verona - ITALY, 
Phone: ++39.045.8098979, FAX: ++39.045.8098982, E-mail:
rocchesso@sci.univr.it} }
\end{center}
}

{\bf Abstract}: {Waveguide Meshes are efficient and versatile models of wave propagation along a multidimensional ideal medium. The choice of the mesh geometry affects both the computational cost and the accuracy of simulations. In this paper, we focus on 2D geometries and use multidimensional sampling theory to compare the square, triangular, and hexagonal meshes in terms of sampling efficiency and dispersion error under conditions of critical sampling. The analysis shows that the triangular geometry exhibits the most desirable tradeoff between accuracy and computational cost. }

{\bf EDICS Category:\\  2-ROOM: Room Acoustics and Acoustic System Modeling}\\
\\
{\bf Permission to publish this abstract separately is granted}.

}

\maketitle

\begin{abstract}{Waveguide Meshes are efficient and versatile models of wave propagation along a multidimensional ideal medium. The choice of the mesh geometry affects both the computational cost and the accuracy of simulations. In this paper, we focus on 2D geometries and use multidimensional sampling theory to compare the square, triangular, and hexagonal meshes in terms of sampling efficiency and dispersion error under conditions of critical sampling. The analysis shows that the triangular geometry exhibits the most desirable tradeoff between accuracy and computational cost. }
\end{abstract}


\section{Introduction}
\PARstart{A}{mong} the techniques for modeling wave propagation in multidimensional media, the {\em Digital Waveguide Meshes} have recently been established~\cite{fonroccim,FontanaRocchessoAcustica,Savioja94,vanduyne932,vanduyne931,vanduyne951,vanduyne96} as intuitive and efficient formulations of finite difference methods~\cite{Strikwerda}.

A Waveguide Mesh (WM) is a discrete--time computational structure that is constructed by  tiling a multidimensional medium into regular elements, each  giving a local description of wave propagation phenomena. This local description is lumped into a waveguide junction~\cite{josMADW,vanduyne931}, which is lossless by construction. Therefore, waveguide meshes are free of numerical losses, even though  lumped passive elements can be explicitely inserted to simulate physical losses. However, wavefronts propagating along  a multidimensional WM are affected by  dispersion error, i.e. different frequencies experience different propagation velocities. Numerical dispersion cannot be completely eliminated, but it can be arbitrarily reduced increasing the density of the elements, and minimized choosing the ``least dispersive geometry''. Moreover, interpolation schemes~\cite{Savioja99acc} or off-line warping techniques~\cite{Savioja99} can be applied to attenuate the effects of numerical dispersion.

In this paper we investigate how the density of waveguide junctions and the sampling frequency affect the signal coming out from the model. As the analysis depends on the geometry of the WM, we focus on 2D media, finding properties for the square, triangular and hexagonal WM (named respectively SWM, TWM and HWM in the following). Such properties allow to calculate the bandwidth of a signal produced by a WM working at a given sampling frequency, once its geometry and the density of its junctions have been determined.   

The paper is structured as follows. Section~\ref{background} provides some background material on waveguide meshes (and their interpretation as finite difference schemes) and multidimensional sampling lattices. Section~\ref{Waveguide_Meshes_as_sampling_schemes} illustrates the spatial sampling efficiency of the three waveguide mesh geometries for signals having circular spatial band shape. In Section~\ref{Signal_time_evolution} we explain how the critical spatial sampling affects the choice of the temporal sampling frequency in non-aliasing conditions. In Section~\ref{perf}, the computational performances of the three geometries are compared under critical sampling conditions.

\section{Background}
\label{background}
\subsection{Digital Waveguides and Waveguide Meshes}
\label{Waveguide_Meshes}
An ideal one-dimensional physical waveguide can be modeled, in discrete time, by means of a couple of parallel delay lines where two wave signals, $s_+$ and $s_-$, travel in opposite directions. Such a structure, based on spatial sampling (with interval $D$) and time sampling (with interval $T$), is called a digital waveguide~\cite{joscmj,SmithInBrandenburg98}. If wave propagation in the physical medium is lossless and non-dispersive with speed $c = D/T$, no error is introduced by the discrete-time simulation as long as $s_+$ and $s_-$ are band limited to a band $B = (2 T)^{-1}$. In this case, the signal $s(x,t)$ along the physical waveguide can be reconstructed with no aliasing error from samples of the wave signals: 
\begin{equation}
s(mD,nT)=s_+(mD,nT) + s_-(mD,nT)\;.
\label{DW}
\end{equation}

$N$ digital waveguide terminations can be connected by means of a  lossless scattering junction~\cite{josMADW,vanduyne931}. Preservation of the total energy in the form of Kirchhoff's node equations leads to the scattering equation
\begin{equation}
  s_{i-} = \frac{2}{N}\,\sum_{k=1}^N s_{k+} - s_{i+} \; \; \; \;  i=1,\ldots,N \; , \label{LSJfirst}
\end{equation}
which allows for calculating the outgoing wave signal to the $i$-th waveguide branch from the $N$ incoming wave signals $s_{1+},\ldots,s_{N+}$, under the assumption of equal wave impedance at the junction for all the waveguides.

A WM, as proposed by Van Duyne and Smith in 1993~\cite{vanduyne931,vanduyne932}, is obtained by connecting unit-length digital waveguide branches by means of lossless scattering junctions. 
For the simulation of uniform and isotropic multidimensional media, a few kinds of geometries, all corresponding to tiling the multidimensional space into regular elements, have been proposed: square~\cite{vanduyne931}, triangular~\cite{fonroccim,vanduyne96}, hexagonal~\cite{vanduyne96} for 2D media such as membranes; rectilinear~\cite{Savioja94} and tetrahedral~\cite{vanduyne96,vanduyne951} for 3D media. The WMs which are considered in this paper (SWM, TWM and HWM) are depicted in figure~\ref{meshes}. Among the proposed geometries, the HWM is peculiar because its $3$--port lossless scattering junctions exhibit two different orientations, and it can be interpreted as two interlaced TWMs (see the junctions marked with $\odot$ in figure~\ref{meshes}.c).
\if F\draft
\begin{figure*}[t]
\center{(a)\hspace{5cm}(b)\hspace{5cm}(c)\vspace{.2cm}\hfill}
\centerline{\hfill
\epsfxsize=4.5cm
{\mbox{{\epsfbox{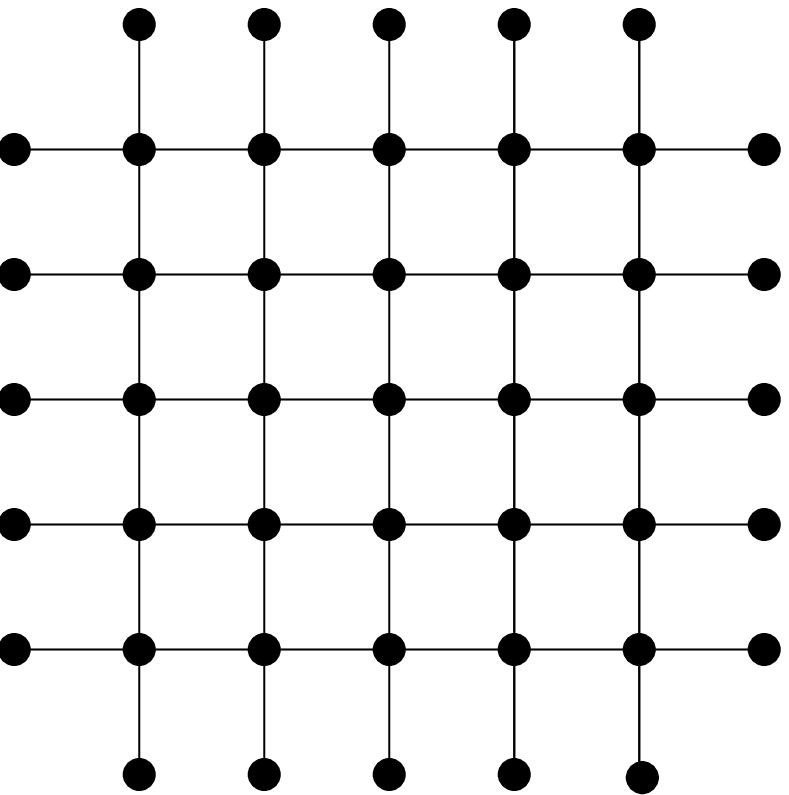}}}}
\hspace{1cm}
\epsfxsize=4.5cm
{\mbox{{\epsfbox{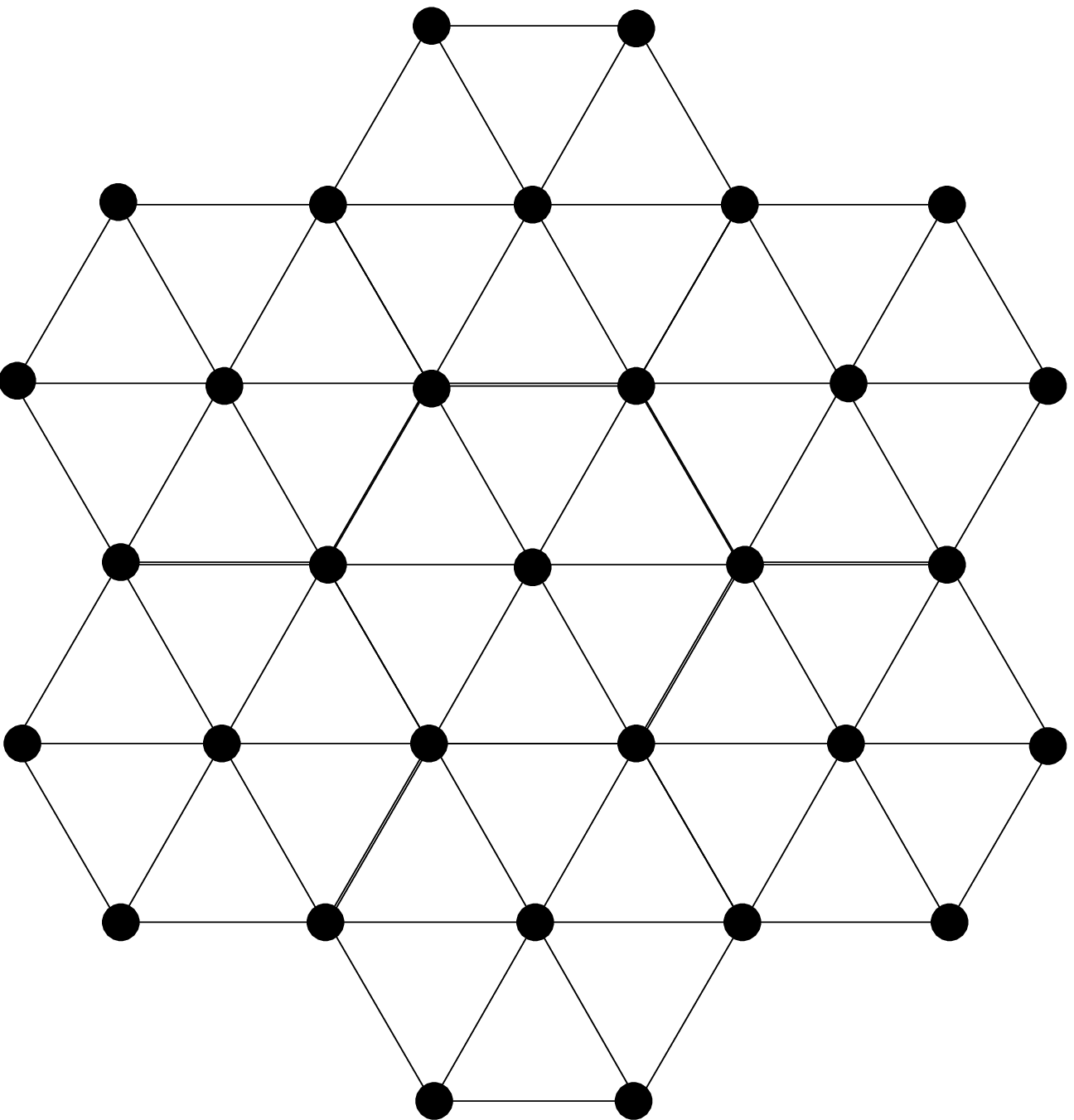}}}}
\hspace{1cm}
\epsfxsize=4.5cm
{\mbox{{\epsfbox{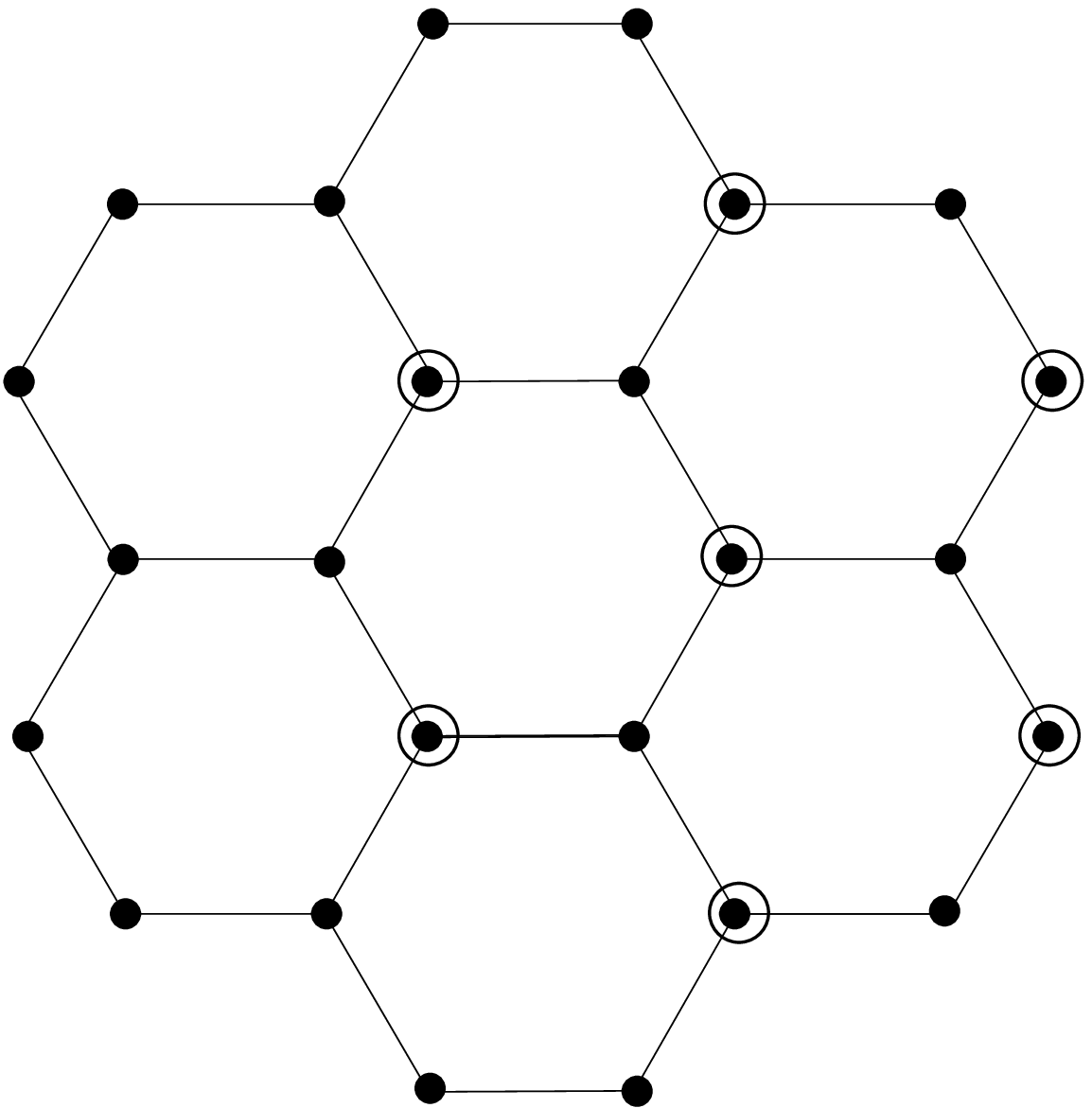}}}}\hfill}
\caption{The SWM (a), the TWM (b) and and the HWM (c). $-$,$/$,$\setminus$ and $|$ are digital waveguides, $\bullet$ are lossless scattering junctions. In (c), seven lossless scattering junctions separated by two digital waveguide branches are marked with $\odot$.}
\label{meshes}
\end{figure*}
\fi

WMs introduce a basic relation, between signal $s({\bf x}_j,nT)$, taken from a junction located at position ${\bf x}_j$ in space, and signals $s({\bf x}_j+{\bf D}_k,nT)\; ,\; k=1,\ldots,N$, taken from the $N$ adjacent junctions connected to it, $D$ meters far from position ${\bf x}_j$: 
\begin{equation}
 s({\bf x}_j,nT+T) + s({\bf x}_j,nT-T) = \frac{2}{N}\sum_{k=1}^N s({\bf x}_j+{\bf D}_k,nT)
\label{WMfund}
\end{equation}
which is obtained from~(\ref{LSJfirst}). 

Equation~(\ref{WMfund}) indicates that each WM behaves like a finite difference scheme~\cite{Strikwerda}, the difference being that the former has a state lumped in the digital waveguides, and the latter has a state lumped in the junctions. For the purpose of the analysis that follows, when the lossless scattering junctions have more than one orientation, as in the HWM, a difference equation such as~(\ref{WMfund}) should be written as many times as there are orientations, each one using a proper set of vectors ${\bf D}_1,\ldots,{\bf D}_N$~\cite{vanduyne96}.

Following the lines of von Neumann stability analysis~\cite{Strikwerda}, equation~(\ref{WMfund}) can be Fourier transformed with spatial variables $x$ and $y$, resulting in
\begin{eqnarray}
S(\xi_x,\xi_y,nT + \alpha_{\mbox{\protect\small g}} T) + S(\xi_x,\xi_y,nT - \alpha_{\mbox{\protect\small g}} T)\;\;\;\;\;\;\;\;\;&& \label{transform}\\
=b_{\mbox{\protect\small g}} S(\xi_x,\xi_y,nT)&& \nonumber
\end{eqnarray}
where $\xi_x$ and $\xi_y$ are spatial frequencies, $\alpha_{\mbox{\protect\small g}}$ takes the value $2$ for the HWM and $1$ for the other meshes, and $b_{\mbox{\protect\small g}}$ is a geometric factor equal to
\begin{eqnarray}
b_{\mbox{\protect\small s}}&=&\cos(2\pi D\xi_x)+\cos(2\pi D\xi_y) \nonumber \\
b_{\mbox{\protect\small t}}&=&\frac{2}{3}\cos\left(2\pi D\xi_x\right)+\frac{2}{3}\cos\left(2\pi D\left[\frac{1}{2}\xi_x+\frac{\sqrt 3}{2}\xi_y\right]\right) \nonumber \\
& & +\frac{2}{3}\cos\left(2\pi D\left[\frac{1}{2}\xi_x-\frac{\sqrt 3}{2}\xi_y\right]\right) \label{geofactor}\\
b_{\mbox{\protect\small h}}&=&\frac{8}{9}\cos\left(2\pi\sqrt{3}D\xi_x\right)  \nonumber \\
& & +\frac{8}{9}\cos\left(2\pi D\left[\frac{\sqrt 3}{2}\xi_x+\frac{3}{2}\xi_y\right]\right)  \nonumber \\
& & +\frac{8}{9}\cos\left(2\pi D\left[\frac{\sqrt 3}{2}\xi_x-\frac{3}{2}\xi_y\right]\right)-\frac{2}{3} \nonumber 
\end{eqnarray}
for the SWM, TWM and HWM, respectively.
\if F\draft
\begin{figure*}[t]
\center{(a)\hspace{5cm}(b)\hspace{5cm}(c)\vspace{.2cm}\hfill}
\centerline{\hfill
\epsfxsize=5.3cm
{\mbox{{\epsfbox{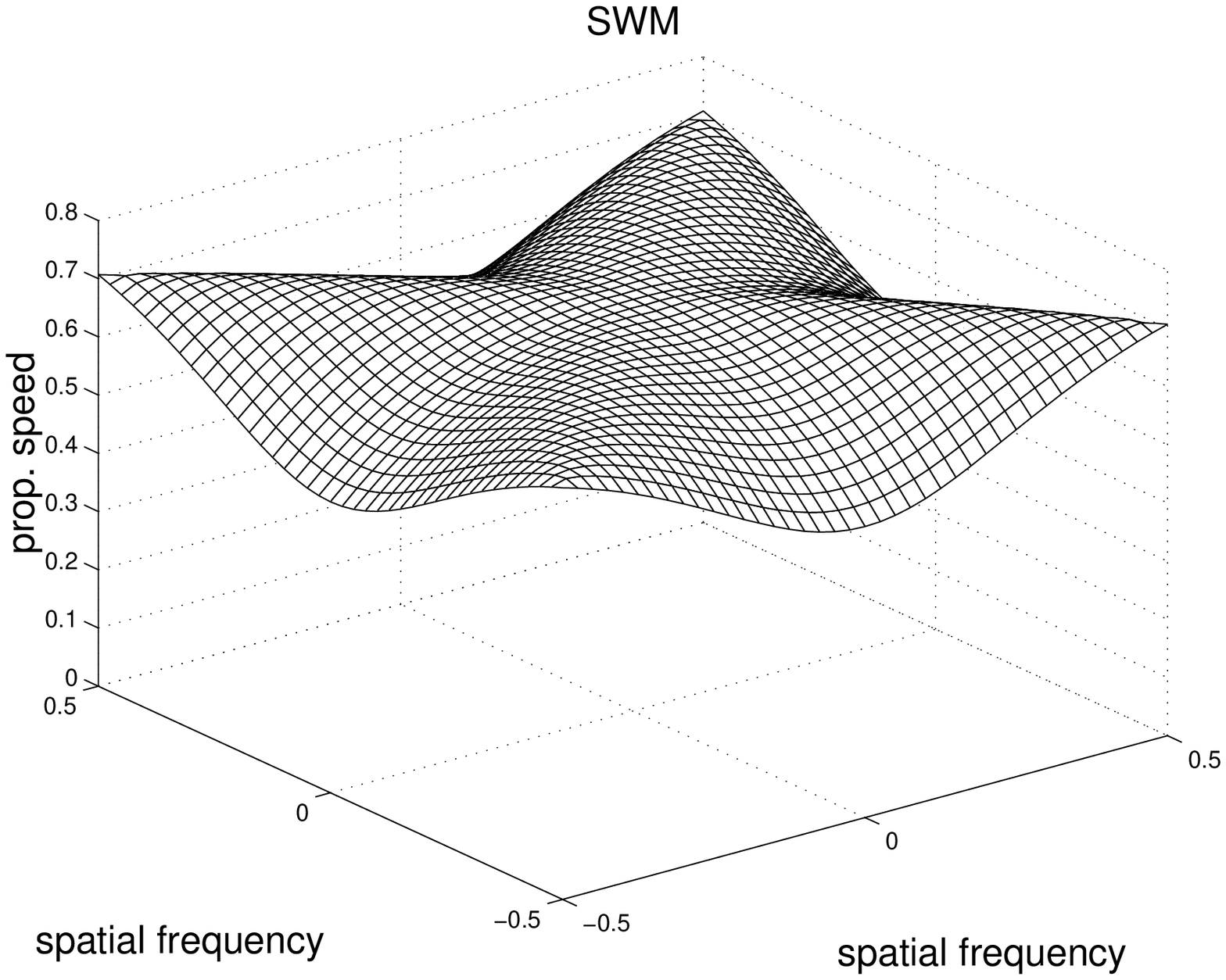}}}}
\hspace{.5cm}
\epsfxsize=5.3cm
{\mbox{{\epsfbox{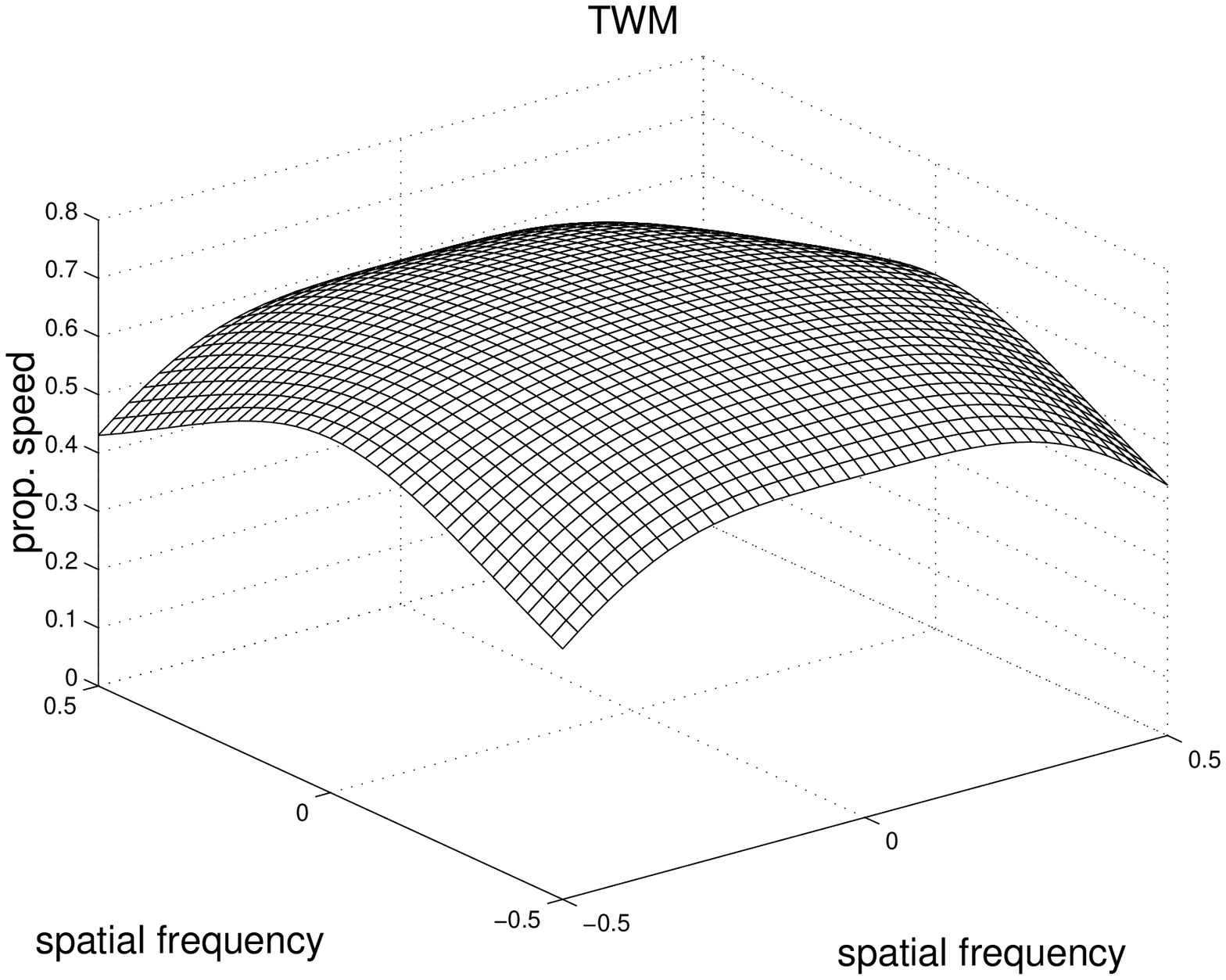}}}}
\hspace{.5cm}
\epsfxsize=5.3cm
{\mbox{{\epsfbox{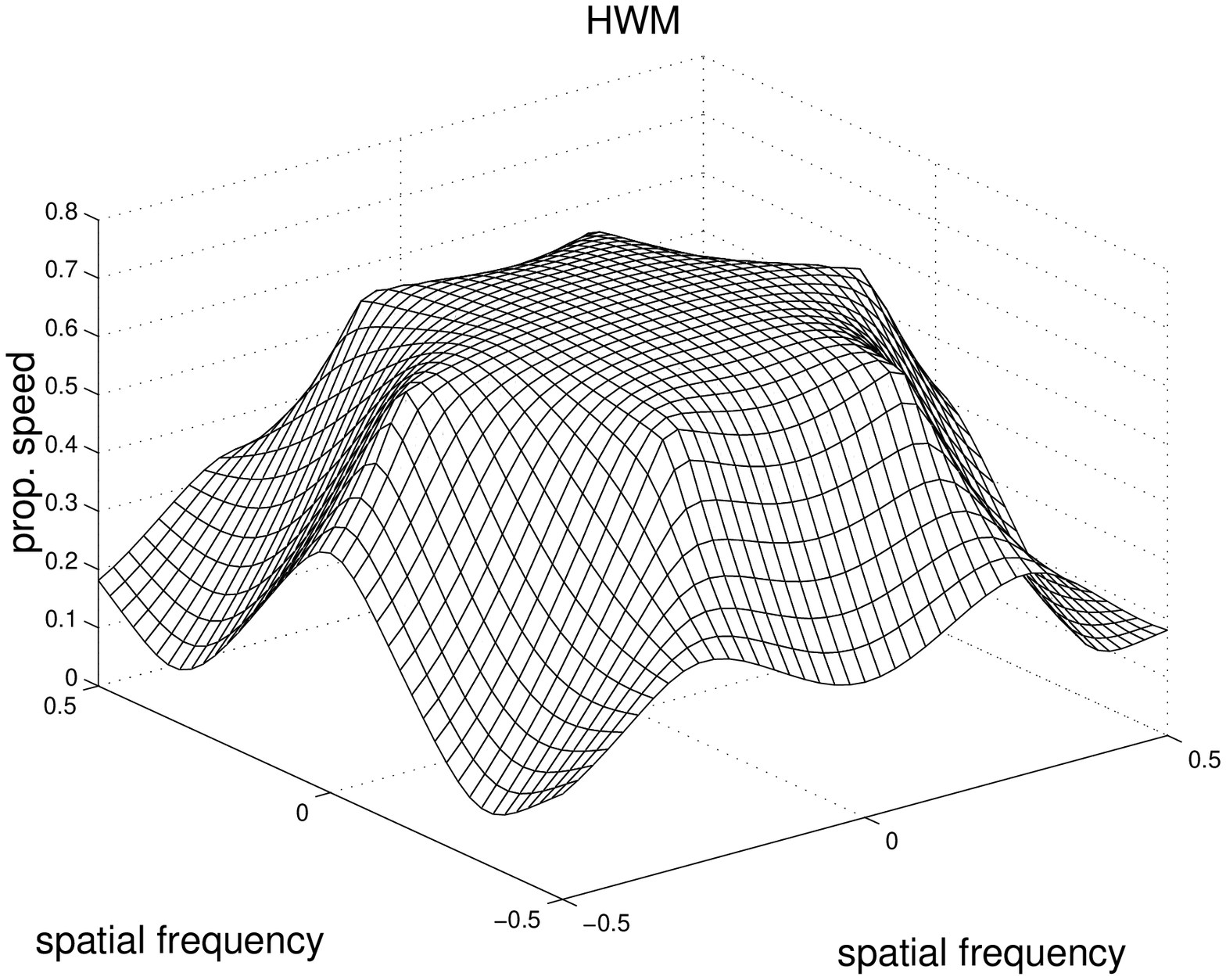}}}}\hfill}
\caption{Propagation speed ratios in the SWM (a), TWM (b) and HWM (c), versus domain given by~(\ref{StartingDomain}).}
\label{DispErr}
\end{figure*}
\fi

Solving equation~(\ref{transform}) as a finite difference equation in the discrete--time variable, the spatial phase shift 
affecting a traveling signal in one time sample is found to be
\begin{equation}
\Comment{\Delta\varphi_{\mbox{g}}(\xi_x,\xi_y)=e^{-j \frac{1}{\alpha_{\mbox {\tiny g}}}\arctan\frac{\sqrt{4-b_{\mbox{\tiny g}}^2}}{b_{\mbox{\tiny g}}}}\: . }
\Delta\varphi_{\mbox{\protect\small g}}(\xi_x,\xi_y)= - \frac{1}{\alpha_{\mbox{\protect\small g}}}\arctan\frac{\sqrt{4-b_{\mbox{\protect\small g}}^2}}{b_{\mbox{\protect\small g}}}\: . 
\label{WMshift} 
\end{equation}

We can compare the propagation speed of a signal traveling along a WM, versus the propagation speed of a signal traveling along an ideal membrane. If we consider membranes where relation $D=cT$ holds, meaning that signals, during a time period $T$, travel for a distance equal to the digital waveguide length, we can calculate the spatial phase shift of these signals, occurring during a time period:
\begin{equation}
\Delta\varphi = - 2\pi D\xi \: ,
\label{idealshift}
\end{equation}
where $\xi=\sqrt{\xi_x^2 + \xi_y^2}$. By comparing~(\ref{WMshift}) and~(\ref{idealshift}), we find the ratio $k_{\mbox{\protect\small g}}$ between the propagation speed of a signal traveling along a WM and along an ideal membrane~\cite{FontanaRocchessoAcustica,Strikwerda,vanduyne932}:
\begin{equation}
k_{\mbox{\protect\small g}}(\xi_x,\xi_y)=\frac{1}{2 \pi \alpha_{\mbox{\protect\small g}} D\xi}\arctan\frac{\sqrt{4-b_{\mbox{\protect\small g}}^2}}{b_{\mbox{\protect\small g}}}  \; . \label{ratios} \\
\end{equation}

Since $k$ is a non constant function of the spatial frequencies, WMs introduce a dispersion error, i.e., the signal traveling along a WM is affected by dispersion of its components. As a starting point, dispersion can be evaluated for spatial frequencies lower than the Nyquist limit\footnote{For the purpose of this paper, the Nyquist limit is defined as half the sample rate, both in time and space. In some previous works~\cite{vanduyne96,Savioja99acc} the simulations were considered valid up to a quarter of the sample rate because in the square mesh the frequency response repeats itself after that limit. However, as it was pointed out in~\cite{vanduyne96}, this is due to the fact that all transfer functions definable at any one junction are functions of $z^{-2}$. This does not imply that the response to a signal having components up to half the sample rate will be aliased. However, for certain applications such as modal analysis of physical membranes, the frequency response at quarter of the sample rate is the definitive limit with the square mesh.}, that is, in the frequency domain:
\begin{equation}
\left\{(\xi_x,\xi_y): |\xi_x|<\frac{1}{2D}\, , \, |\xi_y|<\frac{1}{2D} \right\}\: .
\label{StartingDomain}
\end{equation}
This domain will be refined, according with the considerations to be presented in section~\ref{Waveguide_Meshes_as_sampling_schemes}.

Figure~\ref{DispErr} shows plots of $k_{\mbox{\protect\small s}}$, $k_{\mbox{\protect\small t}}$ and $k_{\mbox{\protect\small h}}$, where $D$ has been set to unity. It can be noticed that the propagation speed decreases for increasing spatial frequencies. In particular, $k_{\mbox{\protect\small s}}$ is maximum when $\xi_x=\xi_y$, and minimum for high values of the spatial frequencies located along the main axes, suggesting that dispersion in the SWM does not affect the diagonal components traveling along it. On the contrary, the HWM exhibits the flattest dispersion error on the region centered around dc. Finally, the TWM seems to have the most uniform behavior of the dispersion error. \Comment{Of course, negative values of the space frequencies have no physical meaning.}

The propagation speed in all the WMs has a maximum at dc:
\begin{equation}
k_{\mbox{\protect\small g}}(0) \define \lim_{\xi\rightarrow 0}{k_{\mbox{\protect\small g}}(\xi)}=\frac{1}{\sqrt 2}\: . \label{DispErrors}  
\end{equation}
It is worth noticing that this value corresponds with the nominal propagation speed of a signal traveling along a finite difference scheme~\cite{Strikwerda}.

\subsection{Sampling Lattices}
The evaluation of the dispersion error does not give a complete description of the constraints holding when an ideal membrane is modeled using WMs. In particular, a method is needed for computing the signal bandwidth a WM is able to process. The theory of sampling lattices~\cite{DudgeonAndMersereau84,vai93}, that we are briefly reviewing in this section, gives the background for characterizing WMs from this viewpoint.

Let us sample a 2D continuous signal $s$ over a domain $L$, subset of ${\cal R}^2$, so defining a discrete signal $s_L({\bf x})\;,\;{\bf x} \in L$. If $L$ can be described by means of a nonsingular matrix ${\bf L}$ such that each element of the domain is a linear combination of the columns of $\bf L$, the coefficients being signed integers:
\begin{equation}
{\bf x} = {\bf L}\;\left[\begin{array}{c}u_1 \\ u_2 \end{array}\right]\:,\:u_1\in{\cal Z}\:,\:u_2\in{\cal Z}\; ,
\label{lattice} 
\end{equation}
then $L$ is called a {\em sampling lattice}, and $\bf L$ is its {\em basis}.

The number of samples per unit area is~\cite{DudgeonAndMersereau84,vai93}:
\begin{equation}
{\cal D}_L=\frac{1}{\det({\bf L})}\; ,
\label{density} 
\end{equation}
as $\bf L$ contains information about the distance between adjacent samples.
 
$S_L$, the Fourier transform of $s_L$, is defined over ${\cal R}^2$ and obtained by periodic imaging of $S$, Fourier transform of $s$. These Fourier images are centered around the elements of the lattice $L^*$, which is described by the basis ${\bf L}^{-T}$, inverse transposed of $\bf L$. Notice that the denser the sampling of $s$ is, the sparser the image centers are. In formulas,
\begin{equation}
{\cal D}_{L^*}=\frac{1}{\det({\bf L}^{-T})}=\det({\bf L})=\frac{1}{{\cal D}_{L}}\; .
\label{densityratios} 
\end{equation}

Let us consider an ideal, unlimited membrane traveled by a spatially band limited signal $s(x,y,t)$, and let us do a spatial sampling of the signal, at a given time. Whenever the spatial sampling defines a sampling lattice $L$, so that $s_L$ is defined, we can calculate $S_L$ using the above results. The multidimensional sampling theorem~\cite{DudgeonAndMersereau84} --- which, in brief, tells that if the Fourier images of $s$ do not intersect one with each other, $s$ can be recovered from $s_L$ with no aliasing error --- indicates whether the chosen sampling scheme induces aliasing.

Equation~(\ref{densityratios}) tells that for each choice of the sampling lattice geometry, the density ${\cal D}_L$ can be increased until the images do not intersect. Conversely, given $s$, there exists a sampling lattice capable to capture all the information needed to recover the original signal using the least density of samples. Clearly, it will exhibit the highest sampling efficiency.

\section{Sampling efficiency of the WMs}
\label{Waveguide_Meshes_as_sampling_schemes}
Sampling efficiency will be calculated for signals having circular spatial band shape centered around the origin of the frequency axes\footnote{This class of signals encompasses the signals obtained when a membrane is excited, in one or several points, by a single shot or by a sequence of shots, using a stick with an approximately round tip.}, with radius equal to $B$. Even if the analysis procedure does not depend on the shape of the spatial domain, the circular band seems to include all the signals occurring in practical applications. Indeed, the procedure is independent of the system evolution --- which in the WMs is controlled by equation~(\ref{LSJfirst}) --- so that it applies to any model discretizing distributed systems, and where signal information can be located over a sampling lattice. This is the case, for example, in finite difference schemes.

\subsection{WMs and Sampling Lattices}
A SWM, having digital waveguides of length $D_{\mbox{\protect\small s}}$, corresponds to a sampling scheme over the lattice ${L}_{\mbox{\protect\small s}}(D_{\mbox{\protect\small s}})$, described by the basis
\begin{equation}
 {\bf L}_{\mbox{\protect\small s}}(D_{\mbox{\protect\small s}}) = \left| 
                      \begin{array}{cc}
                           D_{\mbox{\protect\small s}} & 0 \\
                           0 & D_{\mbox{\protect\small s}}
                      \end{array}
                  \right| \; .
\end{equation}
A TWM, having digital waveguides of length $D_{\mbox{\protect\small t}}$, corresponds to a sampling scheme over the lattice $ {L}_{\mbox{\protect\small t}} ( D_{\mbox{\protect\small t}} ) $, described by the basis
\begin{equation}
 {\bf L}_{\mbox{\protect\small t}}(D_{\mbox{\protect\small t}}) = \left| 
                      \begin{array}{cc}
                           D_{\mbox{\protect\small t}} & \frac{1}{2}D_{\mbox{\protect\small t}} \\
                           0 & \frac{\sqrt 3}{2}D_{\mbox{\protect\small t}}
                      \end{array}
                \right| \; .
\end{equation}

Notice that a triangular scheme is denser than a square scheme made with digital waveguides of the same length. This is confirmed by relation
\begin{equation}
\frac{{\cal D}_{L_{\mbox{\protect\small t}}(D)}}{{\cal D}_{L_{\mbox{\protect\small s}}(D)}}=\frac{\det({\bf L}_{\mbox{\protect\small s}}(D))}{\det({\bf L}_{\mbox{\protect\small t}}(D))}=\frac{D^2}{\frac{\sqrt 3}{2}D^2}=\frac{2}{\sqrt 3} \; .
\end{equation}

The description of an HWM having digital waveguides of length $D_{\mbox{\protect\small h}}$, again can be given in terms of sampling lattices, with some extra care. The HWM is obtained by subtracting a TWM, whose junctions lie on the sampling lattice $L_{\mbox{\protect\small T}}(D_{\mbox{\protect\small h}})$, from a denser TWM having junctions on $L_{\mbox{\protect\small t}}(D_{\mbox{\protect\small h}})$. Hence, the junctions of the HWM lie on
\begin{equation}
L_{\mbox{\protect\small h}}(D_{\mbox{\protect\small h}}) \define L_{\mbox{\protect\small t}}(D_{\mbox{\protect\small h}}) \setminus L_{\mbox{\protect\small T}}(D_{\mbox{\protect\small h}})\; .
\label{hexlattice}
\end{equation}
The basis of $L_{\mbox{\protect\small T}}(D_{\mbox{\protect\small h}})$ is
\begin{equation}
  {\bf L}_{\mbox{\protect\small T}}(D_{\mbox{\protect\small h}}) = \left|
                         \begin{array}{cc}
                           \frac{3}{2}D_{\mbox{\protect\small h}} & 0 \\
                           \frac{\sqrt 3}{2}D_{\mbox{\protect\small h}} & {\sqrt 3}D_{\mbox{\protect\small h}}
                         \end{array}
                 \right|\; . 
\end{equation}
\if F\draft
 \begin{figure}[h]
 \epsfysize=5cm
 \centerline{\hfill\epsfbox{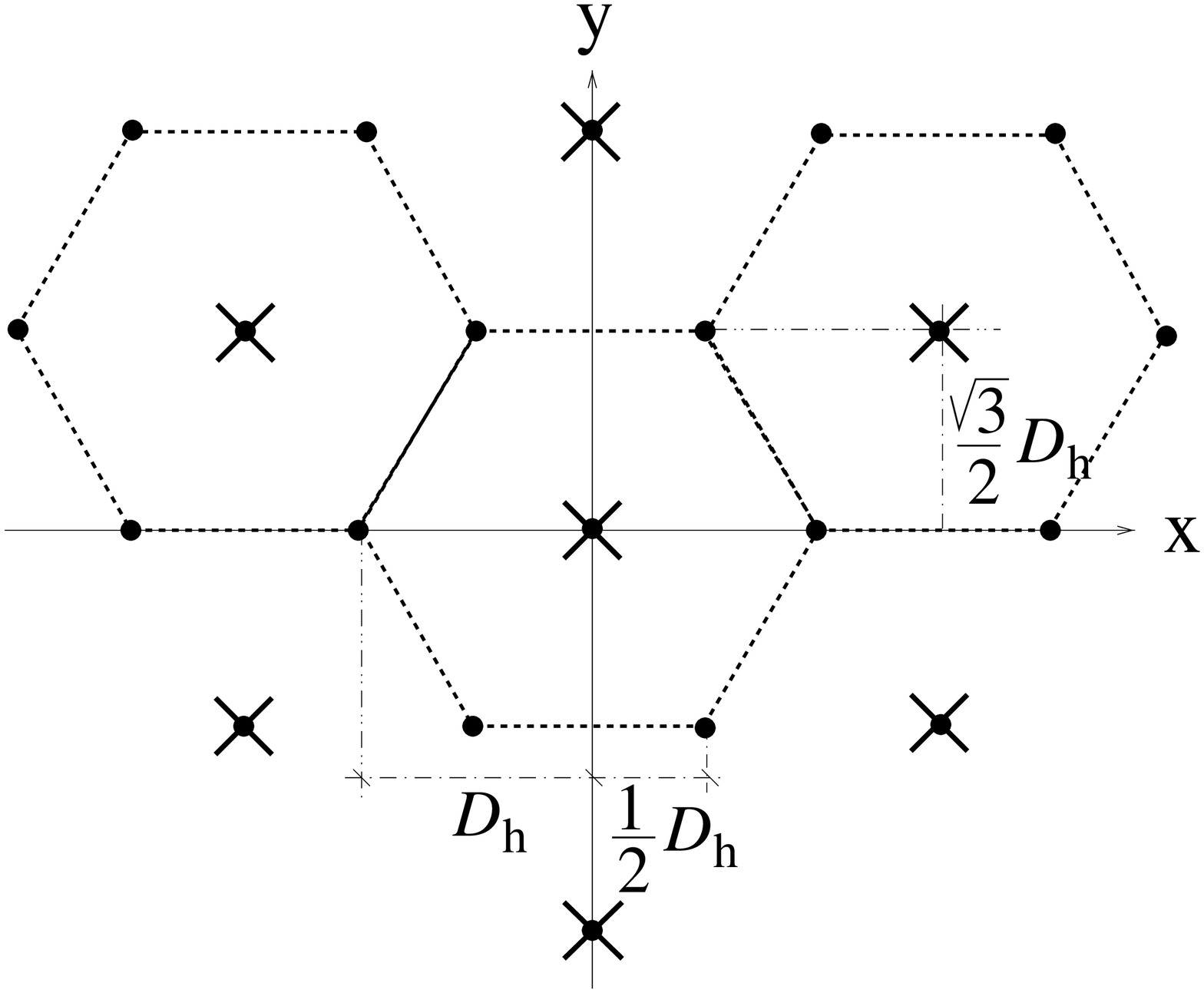}\hfill}
 \caption{The HWM, obtained by subtraction of TWMs. $\times$ are elements belonging to $L_{\mbox{\protect\small T}}(D_{\mbox{\protect\small h}})$, $\bullet$ are elements belonging to $L_{\mbox{\protect\small t}}(D_{\mbox{\protect\small h}})$.}
 \label{HexLatt}
 \end{figure}
\fi

Figure~\ref{HexLatt} shows the HWM over $L_{\mbox{\protect\small h}}(D_{\mbox{\protect\small h}})$, obtained by subtracting the sampling lattice $L_{\mbox{\protect\small T}}(D_{\mbox{\protect\small h}})$ (whose elements are depicted with $\times$) from $L_{\mbox{\protect\small t}}(D_{\mbox{\protect\small h}})$ (elements depicted with $\bullet$).

\subsection{TWM vs SWM}
The image centers of the spectra belonging to signals traveling along a SWM and a TWM are respectively described by the basis matrices of $L^*_{\mbox{\protect\small s}}$ and $L^*_{\mbox{\protect\small t}}$: 
\begin{equation}
  {{\bf L}_{\mbox{\protect\small s}}}(D_{\mbox{\protect\small s}})^{-T} = \left|
                         \begin{array}{cc}
                           1/D_{\mbox{\protect\small s}} & 0 \\
                           0 & 1/D_{\mbox{\protect\small s}}
                         \end{array}
                     \right|\; ,
\end{equation}
and
\begin{equation} 
 {{\bf L}_{\mbox{\protect\small t}}}(D_{\mbox{\protect\small t}})^{-T} = \left|
                      \begin{array}{cc}
                           \frac{1}{D_{\mbox{\protect\small t}}} & 0 \\
                          -\frac{1}{\sqrt 3} \frac{1}{D_{\mbox{\protect\small t}}} & \frac{2}{\sqrt 3} \frac{1}{D_{\mbox{\protect\small t}}}
                      \end{array}
                    \right| \; .
\end{equation}
\if F\draft
 \begin{figure}[h]
 \epsfxsize=220pt
 \centerline{\hfill\epsfbox{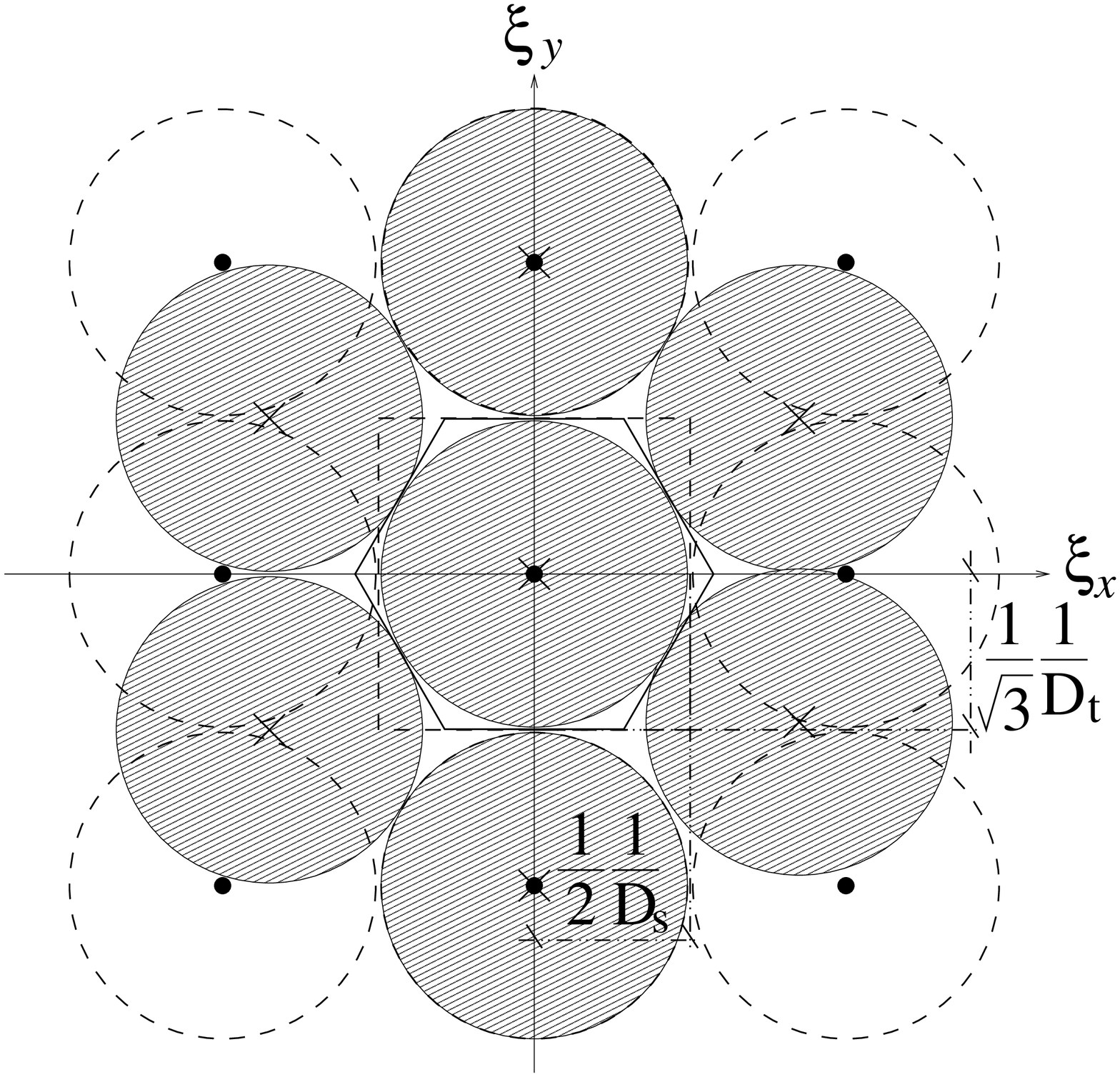}\hfill}
 \caption{Domains of Fourier images in a SWM (empty circles in dashed line) and in a TWM (filled circles), and the correspondent tiling induced by the respective geometries (square in dashed line and hexagon in solid line). The empty and the filled circles have the same radius and touch each other without intersecting, meaning that both the SWM and the TWM critically sample the same signal.}
 \label{EstSqTr}
 \end{figure}
\fi

Figure~\ref{EstSqTr} shows Fourier images for the SWM (empty circles in dashed line) and the TWM (filled circles), located around the origin of the frequency plane. It can be noticed that the square sampling scheme induces a square positioning of the images and, consequently, a square tiling of the frequency plane, as emphasized by the square in dashed line. Similarly, sampling over a TWM results in a triangular positioning of the images, thus producing an hexagonal tiling of the frequency plane, as shown by the hexagon located over the center. Such hexagonal tiling allows, as it appears quite evidently in the figure, to ``pack'' the images better than those coming from a square sampling. 

When both a SWM and a TWM critically sample the same signal (i.e. the filled circles touch each other without intersecting, and the same thing happens for the empty circles), it is easy to derive the following relation between the digital waveguide lengths:
\begin{equation}
\frac{D_{\mbox{\protect\small t}}}{D_{\mbox{\protect\small s}}} = \frac{2}{\sqrt 3}\approx 1.134\; . \label{dtds}
\end{equation}

Under this condition we can relate the sample densities in the two geometries:
\begin{equation}
\frac{{\cal D}_{L_{\mbox{\protect\small t}}(D_{\mbox{\protect\small t}})}}{{\cal D}_{L_{\mbox{\protect\small s}}(D_{\mbox{\protect\small s}})}}=\frac{{\cal D}_{L_{\mbox{\protect\small t}}(\frac{2}{\sqrt 3}D_{\mbox{\protect\small s}})}}{{\cal D}_{L_{\mbox{\protect\small s}}(D_{\mbox{\protect\small s}})}}=\frac{\det({\bf L}_{\mbox{\protect\small s}}(D_{\mbox{\protect\small s}}))}{\det({\bf L}_{\mbox{\protect\small t}}(\frac{2}{\sqrt 3}D_{\mbox{\protect\small s}}))}=\frac{\sqrt 3}{2}\; ,
\label{DensityRatio1} 
\end{equation}
concluding that {\em the TWM exhibits a better sampling efficiency relative to the SWM}. In other words, a signal can be spatially sampled with a triangular geometry using $13.4$\% less samples per unit area.

\subsection{TWM vs HWM}
Said $s_{\mbox{\protect\small h}}$ and $s_{\mbox{\protect\small t}}$ the signals sampled over the junctions of the HWM and the TWM, respectively, and said $s_{\mbox{\protect\small T}}$ the signal sampled over the lattice $L_{\mbox{\protect\small T}}(D_{\mbox{\protect\small h}})$, we can define the zero-padded signals
\begin{equation}
\tilde{s}_{\mbox{\protect\small T}}({\bf x}) = \left\{
   \ba{lcl}
   0           &, & {\bf x} \in L_{\mbox{\protect\small h}}(D_{\mbox{\protect\small h}})  \\
   s_{\mbox{\protect\small T}}({\bf x})&, & {\bf x} \in L_{\mbox{\protect\small T}}(D_{\mbox{\protect\small h}})
   \ea \right.
\end{equation}
and
\begin{equation}
\tilde{s}_{\mbox{\protect\small h}}({\bf x}) = s_{\mbox{\protect\small t}}({\bf x}) - \tilde{s}_{\mbox{\protect\small T}}({\bf x}) =  \left\{
   \ba{lcl}
   s_{\mbox{\protect\small h}}({\bf x})&, & {\bf x} \in L_{\mbox{\protect\small h}}(D_{\mbox{\protect\small h}}) \\
   0           &, & {\bf x} \in L_{\mbox{\protect\small T}}(D_{\mbox{\protect\small h}})
   \ea \right.
\end{equation}

Since we cannot define a Fourier transform over $L_{\mbox{\protect\small h}}(D_{\mbox{\protect\small h}})$, we consider $\tilde{S}_{\mbox{\protect\small h}}$, Fourier transform of $\tilde{s}_{\mbox{\protect\small h}}$ (which is defined over $L_{\mbox{\protect\small t}}(D_{\mbox{\protect\small h}})$), as a description of $s_{\mbox{\protect\small h}}$ in the frequency domain. $\tilde{S}_{\mbox{\protect\small h}}$, according to its definition, can be obtained by subtracting $\tilde{S}_{\mbox{\protect\small T}}$ from $S_{\mbox{\protect\small t}}$:
\begin{equation}
\tilde{S}_{\mbox{\protect\small h}}(\xi_x,\xi_y)= S_{\mbox{\protect\small t}}(\xi_x,\xi_y)-\tilde{S}_{\mbox{\protect\small T}}(\xi_x,\xi_y) \; .
\end{equation}
The result is shown in figure~\ref{EstHex}, where \Comment{the central image of $\tilde{S}_{\mbox{\protect\small h}}$ --- corresponding to the slices, located inside the inner hexagon, of the $6$ grey-filled circles --- is depicted, together with part of its $6$ neighbor images.} some images of $S_{\mbox{\protect\small T}}$ are depicted (filled circles plus circle in dashed line). The elements of ${L_{\mbox{\protect\small T}}}^*(D_{\mbox{\protect\small h}})$ are marked with $\bullet$, while the elements of ${L_{\mbox{\protect\small t}}}^*(D_{\mbox{\protect\small h}})$ are marked with $\times$.
\if F\draft
 \begin{figure}[t]
 \epsfxsize=220pt
 \centerline{\hfill\epsfbox{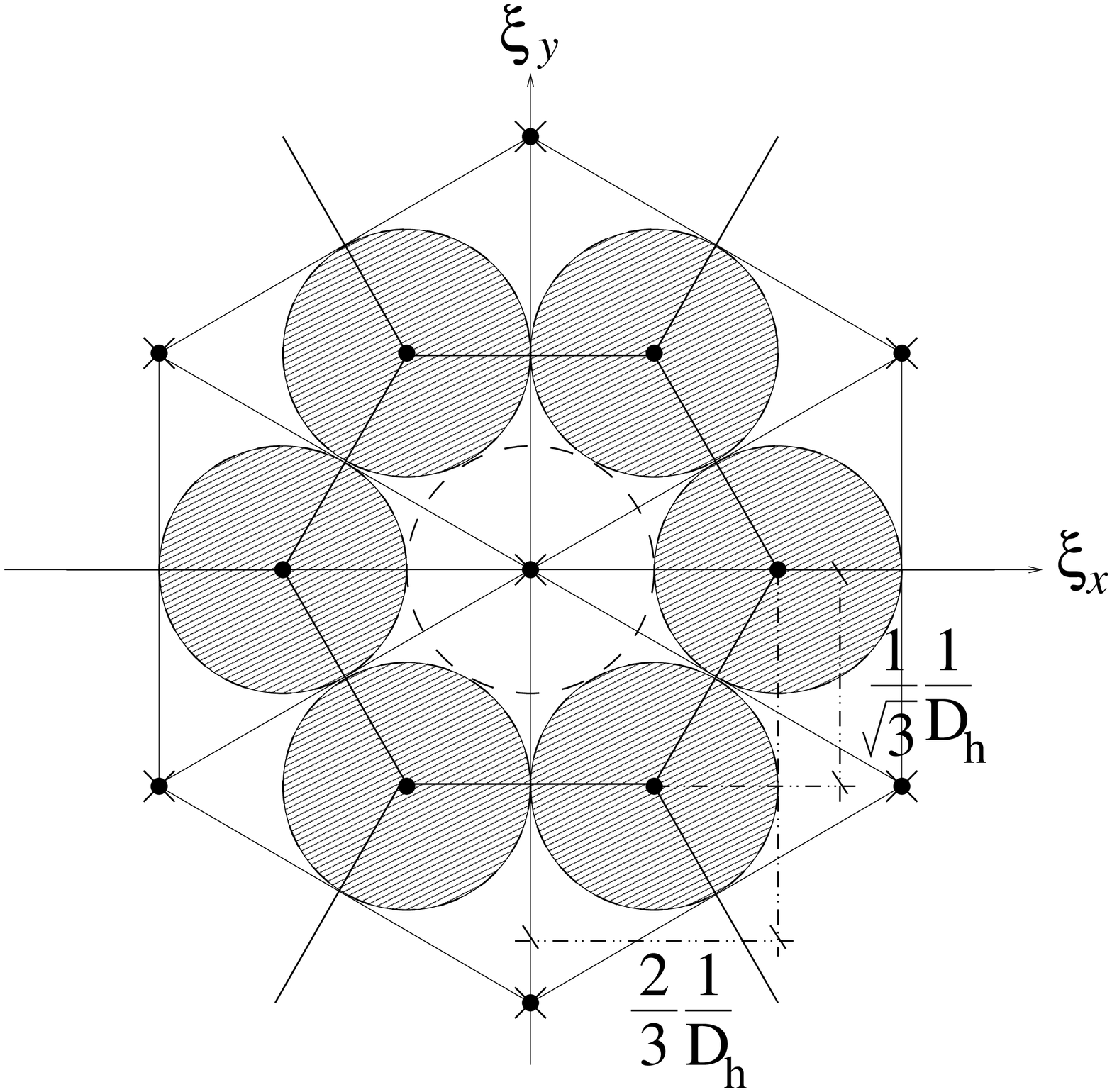}\hfill}
 \caption{Centers ($\times$) of the Fourier images coming from a TWM defined by $L_{\mbox{\protect\small t}}(D_{\mbox{\protect\small h}})$, and centers ($\bullet$) of the Fourier images coming from a sparser TWM defined by $L_{\mbox{\protect\small T}}(D_{\mbox{\protect\small h}})$. Subtraction of $\times$ from $\bullet$ gives the centers of the Fourier images (located on the filled circles) coming from an HWM defined by $L_{\mbox{\protect\small h}}(D_{\mbox{\protect\small h}})$. The respective tiling is given by the triangles. Both the HWM and the sparser TWM critically sample the same signal.}
 \label{EstHex}
 \end{figure}
\fi

The HWM induces an hexagonal positioning of the Fourier images (corresponding to the filled circles). In fact, their centers are elements of a set which, again, can be defined by a subtraction between two sampling lattices, ${L_{\mbox{\protect\small t}}}^*(D_{\mbox{\protect\small h}})$ and ${L_{\mbox{\protect\small T}}}^*(D_{\mbox{\protect\small h}})$, which are reciprocal of ${L_{\mbox{\protect\small t}}}(D_{\mbox{\protect\small h}})$ and ${L_{\mbox{\protect\small T}}}(D_{\mbox{\protect\small h}})$, respectively.
The consequent tiling geometry is triangular, as emphasized by the triangles depicted in figure~\ref{EstHex}.

It can be observed that {\em image intersection does not occur in } $L_{\mbox{\protect\small h}}(D_{\mbox{\protect\small h}})$ {\em if and only if $s$ is sampled without aliasing over } $L_{\mbox{\protect\small T}}(D_{\mbox{\protect\small h}})$. In other words, when $S_{\mbox{\protect\small T}}$ exhibits superposition of its images, the same thing happens for the images of $\tilde{S}_{\mbox{\protect\small h}}$, and vice versa.

Moreover, it must be noticed that equation
\begin{equation}
L_{\mbox{\protect\small T}}(D_{\mbox{\protect\small h}})=L_{\mbox{\protect\small t}}(\sqrt{3}D_{\mbox{\protect\small h}}) \label{dhdh}
\end{equation}
holds between $L_{\mbox{\protect\small T}}$ and $L_{\mbox{\protect\small t}}$, if a rotation is neglected. This relation, together with the considerations about image superposition made just above, allows us to say that an {\mbox{HWM}} does not perform better than a TWM, whose digital waveguides are $\sqrt{3}$ times longer.

Another interesting consideration comes out by noticing that, since
\begin{equation}
\frac{{\cal D}_{L_{\mbox{\protect\small t}}(D_{\mbox{\protect\small h}})}}{{\cal D}_{L_{\mbox{\protect\small T}}(D_{\mbox{\protect\small h}})}}= \frac{\det({\bf L}_{\mbox{\protect\small T}}(D_{\mbox{\protect\small h}}))}{\det({\bf L}_{\mbox{\protect\small t}}(D_{\mbox{\protect\small h}}))} = 3\; ,
\end{equation}
then the TWM defined over the lattice $L_{\mbox{\protect\small T}}(D_{\mbox{\protect\small h}})$ is $3$ times as sparse as the TWM defined over $L_{\mbox{\protect\small t}}(D_{\mbox{\protect\small h}})$, thus (directly from the definition of $L_{\mbox{\protect\small h}}$) $2$ times as sparse as the HWM. Hence, the number of lossless scattering junctions per unit area of the HWM is twice as high as that of the TWM defined over $L_{\mbox{\protect\small T}}(D_{\mbox{\protect\small h}})$, with no benefits on the accuracy of sampling. This relates with the fact that the hexagons tiling the frequency plane (as the inner hexagon depicted in figure~\ref{EstHex}) contain twice the information needed to recover $s$ correctly: the $6$ partial images included in the slices  inside the central hexagon can be composed into $2$ Fourier images.

\section{Signal time evolution}
\label{Signal_time_evolution}
In this section, we discuss the consequences of critical spatial sampling on the temporal sampling frequency.
\if F\draft
\begin{figure*}[t]
\center{(a)\hspace{5.5cm}(b)\hspace{5.5cm}(c)\vspace{.2cm}\hfill}
\centerline{\hfill
\epsfxsize=5cm
{\mbox{{\epsfbox{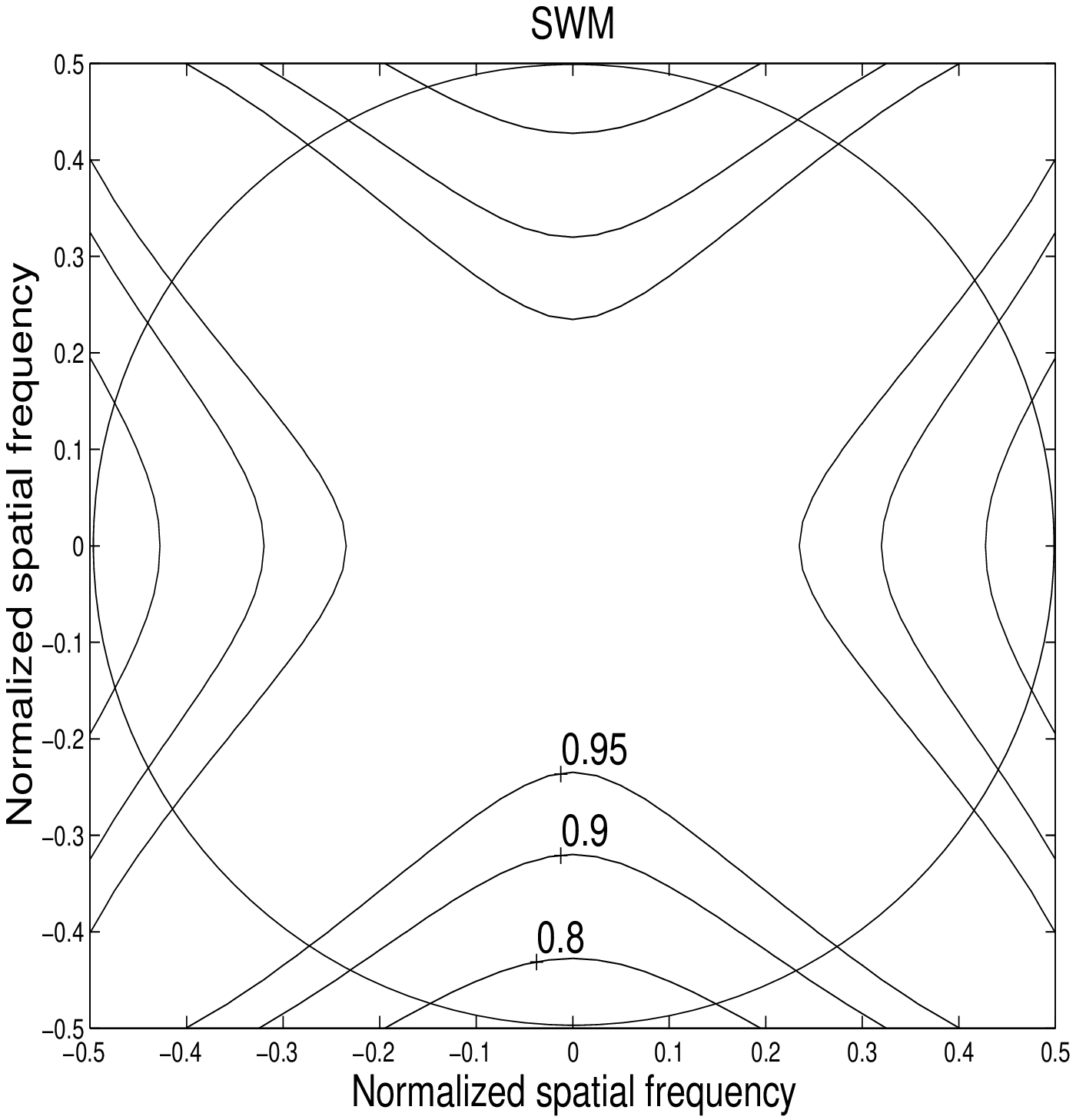}}}}
\hspace{.5cm}
\epsfxsize=5cm
{\mbox{{\epsfbox{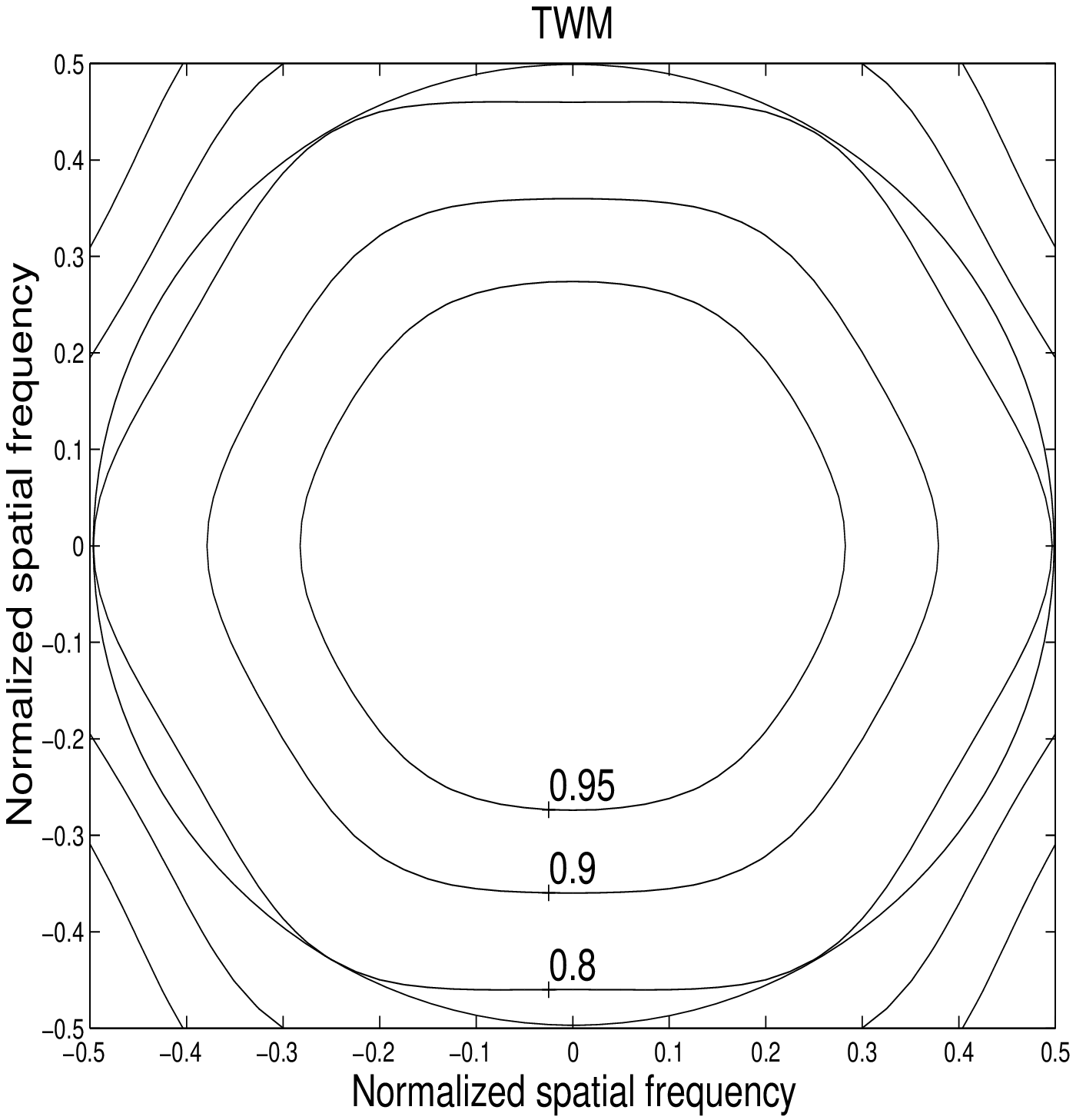}}}}
\hspace{.5cm}
\epsfxsize=5cm
{\mbox{{\epsfbox{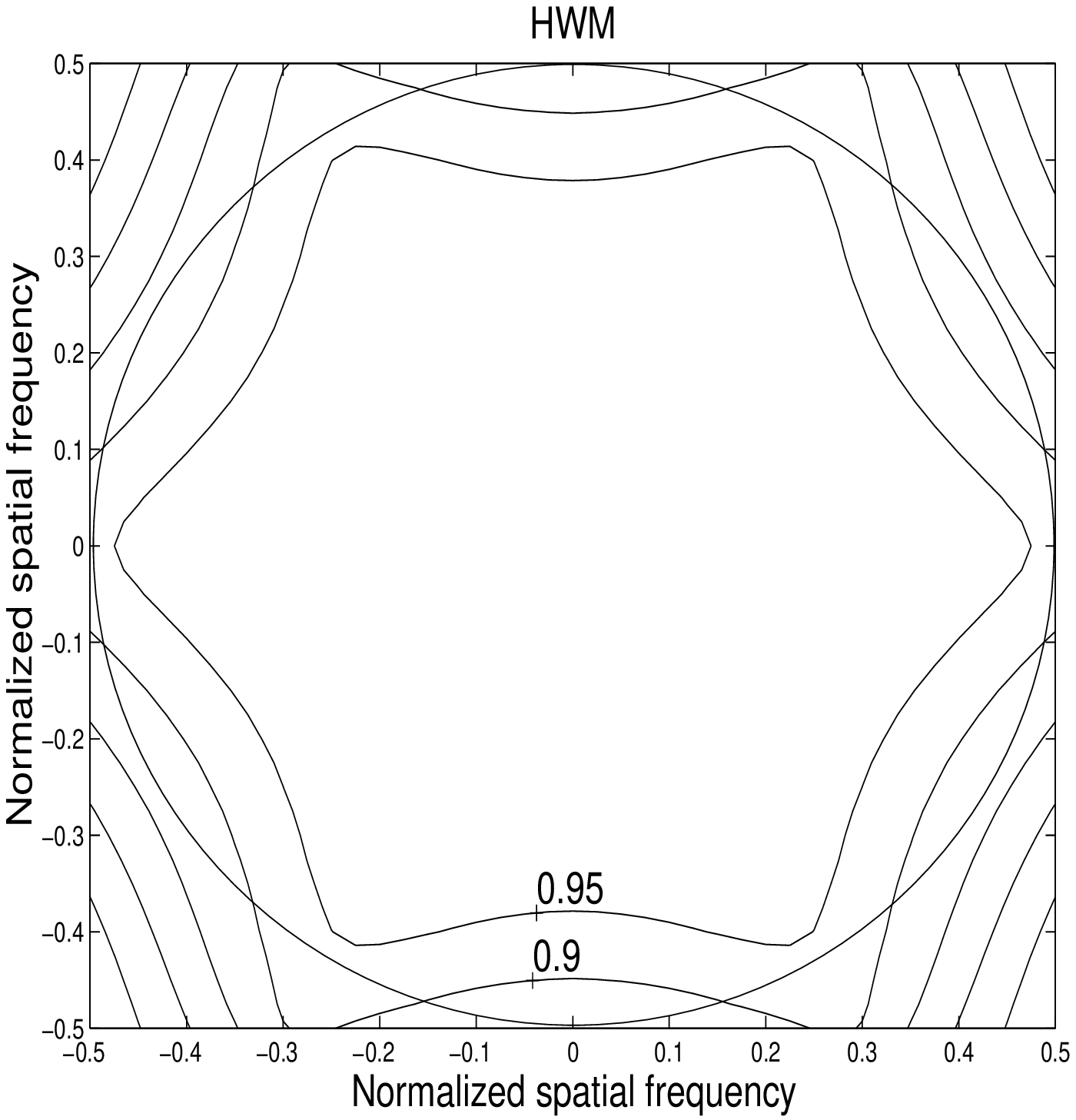}}}}\hfill}
\caption{Contour plots of propagation speed ratios in the SWM (a), TWM (b) and HWM (c), when temporal sampling frequency has been set to $\overline{F}_{s,{\mbox{\protect\small g}}}$. All the WMs process a signal having bandwidth (depicted with the circle) equal to $1/2$.}
\label{CorrErr}
\end{figure*}
\fi

First, let us consider the sampled version of a signal traveling at speed $c$ along an ideal membrane that has been excited by a bandlimited signal $f(x,y)$ having Fourier transform equal to $F(\xi_x, \xi_y)$. It can be shown (see Appendix A) that a time sampling frequency
\begin{equation}
F_s = 2c \max_{\xi:|F(\xi)|\neq 0} {\{\xi\}}
\label{FundTimeGen}
\end{equation}
is required to recover the original signal correctly.

When such a signal has a spatial circular band shape, thus belonging to the class seen in section~\ref{Waveguide_Meshes_as_sampling_schemes}, the relation
\begin{equation}
 \max_{\xi:|F(\xi)|\neq 0} {\{\xi\}}=B
\end{equation}
holds, so that equation~(\ref{FundTimeGen}) simplifies into
\begin{equation}
F_s = 2cB \; .
\label{FundTime}
\end{equation}
Relation~(\ref{FundTime}) has an immediate physical interpretation: waves having wavelength equal to $1/B$, propagating at speed $c$ along a medium, exhibit a temporal frequency equal to $cB$. In order to preserve information in their sampled versions, they must be sampled above twice their frequency.

The value given by equation~(\ref{FundTime}), if used as the sampling frequency of a 2D-resonator  model (realized by means of WMs having critical waveguide lengths), causes inaccurate positioning of the modal frequencies. Assumed that dispersion (see section~\ref{Waveguide_Meshes}) cannot in general be eliminated, as it is a consequence of the finite number of directions a $2$D signal can propagate along a WM\footnote{Savioja and V\"alim\"aki~\cite{Savioja99} have pointed out that a signal, coming out from a TWM, can be frequency warped to reduce the dispersion error. This is made possible by the fact that $k_{\mbox{\protect\small t}}(\xi_x,\xi_y)$ can be approximated with a single--variable function $k_{\mbox{\protect\small t}}(\xi)$.}, the propagation speed can at least be set to its physical  value in  low frequency. This can be done by simply rescaling $F_s$ according to the geometry.

Recalling the procedure leading to equation~(\ref{idealshift}), we can reformulate the spatial phase shift of a signal traveling along the ideal membrane, during a time period $T=1/F_s$. Hence, we find the following expression for the ratio~(\ref{ratios}), reformulated using for each geometry its respective critical waveguide length:
\begin{equation}
\tilde{k}_{\mbox{\protect\small g}}(\xi_x, \xi_y)=\frac{1}{2\pi \alpha_{\mbox{\protect\small g}} D\xi}\arctan\frac{\sqrt{4-\tilde{b}_{\mbox{\protect\small g}}^2}}{\tilde{b}_{\mbox{\protect\small g}}} \; ,
\end{equation}
where $\tilde{b}_{\mbox{\protect\small g}}$ has the structure of $b_{\mbox{\protect\small g}}$ (see equation~(\ref{geofactor})), but the waveguide length $D$ has been replaced by  its critical counterpart $D_{\mbox{\protect\small g}}$.

Recalculation of the limit~(\ref{DispErrors}) gives:
\begin{equation}
\tilde{k}_{\mbox{\protect\small g}}(0)=\frac{D_{\mbox{\protect\small g}}}{D}k_{\mbox{\protect\small g}}(0)=\frac{1}{\sqrt 2}\frac{D_{\mbox{\protect\small g}}}{D}\; ,
\label{corrlimits}
\end{equation}
that is, considering ${D_{\mbox{\protect\small s}}}$ as the reference waveguide length, 
\begin{equation}
\tilde{k}_{\mbox{\protect\small s}}(0)=\frac{1}{\sqrt 2}\; , \; \tilde{k}_{\mbox{\protect\small t}}(0)=\frac{\sqrt 2}{\sqrt 3}\; , \; \tilde{k}_{\mbox{\protect\small h}}(0)=\frac{\sqrt 2}{3}\; .
\label{corrlimitsexp}
\end{equation}

This result shows that signals, under conditions of critical sampling, propagate in the WMs at different speeds, according to the geometry. The dispersion ratios can be set to unity at dc by using the temporal sampling frequencies
\begin{equation}
\overline{F}_{s,{\mbox{\protect\small g}}} = \frac{1}{\tilde{k}_{\mbox{\protect\small g}}(0)}F_s \; .
\label{CorrFreq}
\end{equation}
For the purpose of sec.~\ref{perf}, the dispersion ratio of a critically sampled 2D medium, adjusted to be one at dc, is called $\overline{k}_{\mbox{\protect\small g}}(\xi_x, \xi_y)$.

\section{Performance}
\label{perf}
In order to compare the three geometries under critical sampling conditions, we show in figure~\ref{CorrErr} the contour plots of $\overline{k}_{\mbox{\protect\small g}}$, for a nominal spatial bandwidth $B=1/2$ (corresponding to the circle in the figure).
\if F\draft
\begin{table*}[t]
\caption{Performance of the geometries in terms of memory requirement and computational cost (multiplications are pure bit--shifts in the SWM). Two implementations are considered: as a waveguide mesh and as a finite difference scheme.}
\centering
\begin{tabular}{||c||c|c|c||c|c|c||}                                                                                              \hline 
\Comment{\parbox{1.5in}{\begin{center}{}\end{center}}                                                  & \multicolumn{3}{|p{4.5cm}||}{Waveguide Mesh}  & \multicolumn{3}{|p{4.5cm}||}{Finite Difference Scheme}       \\ \hline \hline}

\parbox{2.5in}{\begin{center}{}\end{center}}                                                  & \multicolumn{3}{|p{3cm}||}{{\begin{center}{Waveguide Mesh}\end{center}}}  & \multicolumn{3}{|p{3cm}||}{{\begin{center}{Finite Difference Scheme}\end{center}}}       \\ \hline \hline

\parbox{2.5in}{\begin{center}{}\end{center}}                                                  & Sq.     & Tr.         & Hex.  & Sq.     & Tr.         & Hex.       \\ \hline \hline
\parbox{2.5in}{\begin{center}{\em additions per junction}\end{center}}                        &   $7$  &   $11$      &  $5$  & 4 & 6 & 3     \\  \hline
\parbox{2.5in}{\begin{center}{\em multiplications per junction}\end{center}}                  &   $1$  &   $1$      &  $1$  & 1 & 1 & 1     \\  \hline
\parbox{2.5in}{\begin{center}{\em memory locations per junction}\end{center}}                 &   $4$  &   $6$       &  $3$ & 2 & 2 & 2      \\  \hline
\parbox{2.5in}{\begin{center}{\em density of junctions}\end{center}}                          &   $1$  &   $0.866$   &  $1.732$ & $1$  &   $0.866$   &  $1.732$  \\ \hline
\parbox{2.5in}{\begin{center} {\em density of memory locations}\end{center}}                  &   $4$  &   $5.2$     &  $5.2$  & $2$ & $1.732$ & $3.464$    \\ \hline
\parbox{2.5in}{\begin{center}{\em sample rate}\end{center}}                                   &   $1$  &   $0.866$   &  $1.5$ & $1$  &   $0.866$   &  $1.5$     \\ \hline
\parbox{2.5in}{\begin{center}{\em additions per unit time and space}\end{center}}             &   $7$  &   $8.25$    &  $12.990$ & $4$ & $4.5$ & $7.794$ \\ \hline
\parbox{2.5in}{\begin{center}{\em multiplications per unit time and space}\end{center}}       &   $1$  &   $0.75$    &  $2.598$  & $1$  &   $0.75$    &  $2.598$\\ \hline
\end{tabular}
\label{performance}
\end{table*}
\fi

It can be noted that the behaviors of the SWM and the TWM are not dramatically different in terms of average dispersion error: dispersion stands quite below $20\%$ for most spatial frequencies in both the geometries. However, the TWM exhibits a more uniform behavior, and this uniformity can be exploited using frequency warping~\cite{Savioja99acc}. Conversely, dispersion in the HWM stays below\Comment{Sugg.4 of rew.2, 3/2000} $10\%$ almost everywhere, thus indicating that this geometry has the most uniform propagation speed under these test conditions.

The computational cost and the memory requirement in the different geometries can be calculated under the same conditions. They are based on equation~(\ref{LSJfirst}) which shows that each $N$--port lossless scattering junction requires $2N$ operations to compute the wave signals coming out from the junctions, to be stored into $N$ locations belonging to the adjacent digital waveguides. These operations amount to $2N\!-\!1$ additions, and $1$ multiplication (which can be replaced by a bit shift in fixed--point implementations of the SWM). 

Table~\ref{performance} summarizes the performance, when the reference sampling rate of the square mesh ${\overline F}_{s,{\mbox{\protect\small s}}}$ has been set to a nominal value equal to unity. Numbers are given for two implementations: as a waveguide mesh (with memory in the waveguide branches), and as a finite difference scheme (with memory in the junctions).

The numbers of operations and memory locations per junction for the WM follow directly from equation~(\ref{LSJfirst}). The numbers of operations and memory locations per junction for the finite difference scheme follow directly from equation~(\ref{WMfund}). The densities of junctions result from equations~(\ref{dtds}) and~(\ref{dhdh}). The densities of locations are obtained by multiplying the third row times the fourth row. The sample rates are a consequence of equation~(\ref{CorrFreq}).

The results of table~\ref{performance} show that the triangular finite difference scheme uses the least quantity of memory. Among WM implementations, the SWM is the most efficient in terms of memory occupation.

Finally, the number of additions (multiplications) per unit time and space results by multiplication of the number of additions (multiplications) per junction, density of junctions, and sample rate. Once again, the SWM has the least density of operations when a WM model is required. Conversely, the square and the triangular geometries have about the same computational requirements in a finite difference implementation.

\subsection{Numerical example}
Let us  design a WM, capable of modeling in real--time an ideal $2$D round resonator, of radius $r=0.1$ m, where waves propagate at a speed $c$ equal to $130$ m$/$s. Let the signal contain information up to a frequency $f=10$ kHz.

The spatial bandwidth $B$ of the signal traveling along the resonator is:
\begin{equation}
B=\frac{f}{c}=\frac{10000\; {\mbox{Hz}}}{130\; {\mbox{m}}/{\mbox{s}}}=76.923\; {\mbox m}^{-1}\; ,
\end{equation}
and consequently the critical waveguide lengths required to model the signal, in the respective geometries, are:
\beqa
D_{\mbox{\protect\small s}} & = & \frac{1}{2B}=6.5\; \mbox{mm} \; ; \nonumber \\
D_{\mbox{\protect\small t}} & = & \frac{1}{\sqrt{3}B}=7.5\; \mbox{mm} \; ;  \\
D_{\mbox{\protect\small h}} & = & \frac{1}{3B}=4.3\; \mbox{mm}\; . \nonumber
\eeqa
This means that the numbers of junctions --- $N_{\mbox{\protect\small s}}$, $N_{\mbox{\protect\small t}}$ and $N_{\mbox{\protect\small h}}$, respectively --- needed in the three models are\footnote{The values result by calculating the number of small areas, each one being associated with its own lossless scattering junction, which tessellate the resonator.}
\beqa
N_{\mbox{\protect\small s}} & \approx & \frac{\pi r^2}{D_{\mbox{\protect\small s}}^2}=744 \; ; \nonumber \\
N_{\mbox{\protect\small t}} & \approx & \frac{\pi r^2}{\frac{\sqrt 3}{2}D_{\mbox{\protect\small t}}^2}=645 \; ; \\
N_{\mbox{\protect\small h}} & \approx\ & \frac{2}{3}\frac{\pi r^2}{\frac{\sqrt 3}{2}D_{\mbox{\protect\small h}}^2}=1308 \; . \nonumber
\eeqa

The time sampling frequencies needed to have $\overline{k}_{\mbox{\protect\small g}}(0)=1$ are, from equation~(\ref{CorrFreq}):
\beqa
{\overline F}_{s,{\mbox{\protect\small s}}} & = & 2\sqrt{2}f=28.285\; \mbox{kHz}\; ; \; \nonumber \\
{\overline F}_{s,{\mbox{\protect\small t}}} & = & 2\sqrt{\frac{3}{2}}f= 24.495\; \mbox{kHz}\; ; \; \\
{\overline F}_{s,{\mbox{\protect\small h}}} & = & 2\frac{3}{\sqrt 2}f= 42.427\; \mbox{kHz}\; . \nonumber
\eeqa

\section{Conclusion}
A novel wave propagation model has recently been introduced for the simulation of isotropic multidimensional media. It makes use of structures called waveguide meshes. Even if most of the waveguide mesh properties have already been understood, there was lack of literature about their performance from a signal--sampling viewpoint.

In this paper, some properties of the most common $2$D waveguide meshes --- square, triangular and hexagonal --- related with the bandwidth of the signal traveling on them, have been inspected. In particular, it has been shown that the triangular waveguide mesh is capable of processing a larger bandwidth than the square or hexagonal waveguide meshes having the same digital waveguide lengths.

Furthermore, when processing signals of the same bandwidth, the triangular waveguide mesh does not exhibit a computational and memory load much larger than the most computationally--efficient waveguide mesh --- the square mesh --- such that it can be considered, for the uniformity of its dispersion error, a good choice both in terms of simulation errors and computational cost. 

The analysis presented in this article may be extended to $3$D media, comparing the three geometries which most directly correspond to the waveguide meshes here reviewed: $3$D rectilinear, dodecahedral, and tethraedral.

\section{Acknowledgment}
We would like to thank Lauri Savioja and Vesa {V\"alim\"aki} for many insightful discussions. We are also grateful to the anonymous reviewers for their constructive criticism.

\section*{Appendix A: Sampling of a signal traveling along an ideal membrane}
\label{appendice}
An ideal membrane establishes a relation between the spatial Fourier transforms of the signal $s$ traveling on it, taken in correspondence of two times, $\tilde t$ and $t$:
\begin{equation}
S(\xi_x, \xi_y, t)= e^{j2\pi c(t-{\tilde t})\xi}\, S(\xi_x, \xi_y, {\tilde t})\; .
\label{app1}
\end{equation}
This relation has the following interpretation: each spatial component of the signal, traveling along a distance equal to $c(t-{\tilde t})$ during time $({\tilde t},t)$, has no magnitude variation, and has phase variation equal to $-c(t-{\tilde t})\xi$~\cite{vanduyne932}.

Let us excite, at time ${\tilde t}$, the ideal membrane with a spatially band limited signal $f(x,y)$ having Fourier transform  $F(\xi_x,\xi_y)$. The application of~(\ref{app1}) gives the spatial Fourier transform of the signal on the membrane:
\begin{eqnarray}
S(\xi_x, \xi_y, t)   &=&  \left\{
\begin{array}{lcl}
          0                   &\!\! , & t<{\tilde t} \\
          e^{j2\pi c(t-{\tilde t})\xi}\, F(\xi_x, \xi_y) &\!\! , & t\geq {\tilde t}
\end{array}
\right.
\label{spatialtrans}
\end{eqnarray}

From this relation, we can calculate the critical temporal sampling frequency required to sample, without loss of information, a signal traveling on an ideal membrane. In fact, the Fourier transform of~(\ref{spatialtrans}),
\begin{equation}
S(\xi_x,\xi_y,f)=\int_{-\infty}^{+\infty}S(\xi_x,\xi_y,t)e^{-j2\pi ft}\,dt \; ,
\end{equation}
has a magnitude equal to 
\begin{eqnarray}
\lefteqn{|S(\xi_x, \xi_y, f)|} \nonumber \\
&&=\left|\int_{-\infty}^{+\infty} {S(\xi_x, \xi_y, t) e^{-j2\pi ft}\, dt}\right| \nonumber \\
&&=\left|\int_{\tilde t}^{+\infty}e^{j2\pi c(t-{\tilde t})\xi}\, F(\xi_x, \xi_y)e^{-j2\pi ft}\, dt\right| \\
&&=\left|\int_{\tilde t}^{+\infty}e^{-j2\pi c{\tilde t}\xi}\, F(\xi_x, \xi_y)e^{-j2\pi t\left\{ f-c\xi\right\}}\, dt\right| \nonumber \\
&&=\left|F(\xi_x, \xi_y)\right|\left|\int_{\tilde t}^{+\infty}e^{-j2\pi t\left\{ f-c\xi\right\}}\, dt\right|\; .\nonumber 
\end{eqnarray}
In particular, when ${\tilde t}\rightarrow -\infty\,$:
\begin{eqnarray}
\lefteqn{|S(\xi_x, \xi_y, f)|} \nonumber \\
&&=\left|F(\xi_x, \xi_y)\right|\lim_{{\tilde t}\rightarrow -\infty}\left|\int_{\tilde t}^{+\infty}e^{-j2\pi t\left\{ f-c\xi\right\}}\, dt\right| \\
&&=\left|F(\xi_x, \xi_y)\right|\left|\int_{-\infty}^{+\infty}e^{-j2\pi t\left\{ f-c\xi\right\}}\, dt\right| \nonumber \\
&&=\left|F(\xi_x, \xi_y)\right| \,\, \delta\!\left( f-c\xi\right) \; .\nonumber 
\end{eqnarray}
This means that some spectral power exists for any temporal frequency (equal to $c\xi$) associated to a spatial component of frequencies $(\xi_x,\xi_y)$ excited in the membrane. Invoking the sampling theorem~\cite{OppenheimAndSchafer}, the critical time sampling frequency $F_s$ results to be equal to~(\ref{FundTimeGen}).

\bibliography{general}

\Comment{
\if T\draft
\newpage
\input figures.tex

\newpage
\input tables.tex

\newpage
\listoffigures

\newpage
\listoftables
\fi
}

\newpage

\begin{biography}{Federico Fontana} received the {\it Laurea in Ingegneria Elettronica} degree from the University of Padova in 1996. After a period spent in Spain at UH S.A. working in the field of acoustic \& vibration, since 1998 he has been collaborating with the {\it Centro di Sonologia Computazionale} (CSC) at the University of Padova, working on digital equalization of small rooms and sound synthesis by physical modeling. He has been a consultant for Generalmusic and for the DPG Audio \& Automotive Division of STMicroelectronics. He is currently a  Ph.D. student at the University of Verona. His main interests are in sound processing, physical modeling of multidimensional resonators, and jazz drumming as an amateur player.
\end{biography}
\begin{biography}{Davide Rocchesso} received the {\it Laurea in Ingegneria Elettronica} degree from the University of Padova in 1992, and the Ph.D. degree from the same university in 1996. His Ph.D. research involved the design of structures and algorithms based on feedback delay networks for sound processing applications. 
In 1994 and 1995 he was a visiting scholar at the Center for Computer Research in Music and Acoustics (CCRMA) at Stanford University. Since 1991 he has been collaborating with the {\it Centro di Sonologia Computazionale} (CSC) at the University of Padova as a researcher and a live-electronic designer. Since march 1998 he has been with the {\it Dipartimento Scientifico e Tecnologico} at the University of Verona, as an Assistant Professor. His main interests are in sound processing, physical modeling, sound reverberation and spatialization, multimedia systems. His home page on the web is http://www.sci.univr.it/\~{}rocchess.
\end{biography}

\end{document}